\documentclass{amsart}
\usepackage{epsfig}

\newtheorem{thm}{Theorem}[section]
\newtheorem{lemma}[thm]{Lemma}
\newtheorem{proposition}[thm]{Proposition}
\newtheorem{prop}[thm]{Proposition}
\newtheorem{conjecture}[thm]{Conjecture}

\newtheorem*{remark}{Remark}

\newcommand{\BP}[1]{BP\left< #1 \right>}
\DeclareMathOperator{\Ext}{Ext}
\newcommand{\tmf}{\mathit{tmf}}
\DeclareMathOperator{\eo}{eo}
\newcommand{\ko}{\mathit{ko}}
\DeclareMathOperator{\Z}{\!\mathbb Z\!}
\DeclareMathOperator{\sA}{\mathcal A}
\DeclareMathOperator{\zZ}{\!\mathbb Z}

\title{The $3$-local $\tmf$ homology of $B\Sigma_3$}
\author{Michael A.~Hill}

\begin{document}
\bibliographystyle{amsplain}
\begin{abstract}
In this paper, we introduce a Hopf algebra, developed by the author
and Andr\'e Henriques, which is usable in the computation of the
$\tmf$ homology of a space. As an application, we compute the
$\tmf$ homology of $B\Sigma_3$ in a manner analogous to
Mahowald's computation of the $\ko$ homology $\mathbb RP^{\infty}$
in \cite{OpsSq4}.
\end{abstract}
\maketitle

\section{Introduction}

In this paper we compute the 3-local $\tmf$ homology and $\tmf$ Tate
cohomology of the symmetric group $\Sigma_3$.  This computation is
motivated as follows. Mahowald's computation of $\ko_\ast(\mathbb RP^{\infty})$ has proved useful in a variety of contexts. In particular, Mahowald used $\ko_\ast(\mathbb RP^n)$ and $\ko_\ast(\mathbb RP^{\infty}/\mathbb RP^k)$ to get information about $v_1$ metastable homotopy theory in the $EHP$ sequence \cite{ImJ}. Mahowald has also used $\ko_\ast(\mathbb RP^{\infty})$ to detect elements in his $\eta_j$ family \cite{Etaj}. At the prime $3$, the role of the spectrum $\ko$ is most naturally played by the spectrum $\tmf$. To generalize these results of Mahowald's, the initial piece of data needed is the $\tmf$ homology of $B\Sigma_3$. Both of the aforementioned results should be generalizable starting from this point.

A theorem of Arone and Mahowald shows that $v_n$ periodic information is captured by the first $p^n$ stages of the Goodwillie tower \cite{AroneMahowald}. This recasts Mahowald's result from \cite{ImJ} into a more readily generalizable form. To get $v_2$ periodic information at the prime $3$, the initial data needed comes in part from $QS^0$ and $Q({B\Sigma_3}_{k}^{\infty})$, where ${B\Sigma_3}_{k}^{\infty}$ is a particular Thom spectrum of $B\Sigma_3$. Just as Mahowald uses knowledge of the $\ko$ homology of stunted projective spaces to reduce the questions involved to ones of $J$ homology, we hope that a similar analysis, using Behrens' $Q(2)$, spectrum will allow an analysis of the $v_2$ primary Goodwillie tower at $3$ \cite{ModularK(2)}.

Minami shows that the odd primary $\eta_j$ family will be detectable in the Hurwicz image of the $\tmf$ homology of the $n$-skeleton of $B\Sigma_3$ for appropriate choices of $n$ \cite{MinamiEtaj}. While determining the full Hurwicz image is a trickier task, understanding the groups and simple $\tmf$ operations on them could help determine if the conjectural $\eta_j$ elements actually survive at the prime $3$.

\subsection{Organization of Paper}
In \S \ref{HopfAlg}, we introduce the main computational Hopf algebra $\sA$, $\Ext$ over which is the Adams $E_2$ term for computing $\tmf$ homology. In \S \ref{Review}, we review Mahowald's computation of the $\ko$ homology of $\mathbb RP^{\infty}$, presenting it in a manner which can be most readily generalized. In \S \ref{transfer}, we carry out one of the computational steps analogous to Mahowald's, computing the $\tmf$ homology of the cofiber of the transfer map, and in \S \ref{BSig3}, we complete the computation of $\tmf_\ast(B\Sigma_3)$. Rounding out the computations, in \S \ref{Truncated}, we compute the $\tmf$ homology of the finite skeleta of $R_3$, giving additional results about that of the finite skeleta of $B\Sigma_3$.

The last two sections present conjectures as to further results. A computation of the homotopy of the $\Sigma_3$ Tate spectrum for $\tmf$ is presented together with a non-splitting conjecture in \S \ref{tate}. Conjectural generalizations to primes bigger than $3$ are presented in \S \ref{eop-1}, together with the implications to the $\eo_{p-1}$ homology of $B\Sigma_p$.

\subsection{Conventions and Notation}
We restrict attention to the prime $3$ and assume that all spaces and spectra are $3$-completed except in \S \ref{Review}. For ease of readability, let $H$ be $H\Z/3$. If $X$ is a space or spectrum, let $X^{[n]}$ denote its $n$-skeleton. 

Finally, we need some $\tmf$ specific notation. Let $I$ denote the ideal of the Adams $E_2$ term for $\tmf_\ast$ generated by $v_0$, $c_4$ and $c_6$. Let $\bar{I}$ denote the ideal of $\tmf_\ast$ generated by $3$, $c_4$, $c_6$, and their $\Delta$ and $\Delta^2$ translates. The ideal $I$ converges to the ideal $\bar{I}$, and $I$ is the annihilator ideal of the elements $\alpha$ and $\beta$. For brevity, the reader is asked to always assume the relation $I(\alpha, \beta)$ in all Adams $E_2$ terms, unless explicitly stated otherwise. Moreover, we assume that the relations $c_4^3-c_6^2=27\Delta$ always holds and will not be explicitly stated.

\section{Fundamental Hopf Algebra}\label{HopfAlg}

Our basic tool of computation will be a variant of the Adams
spectral sequences based on ordinary cohomology.
Since $H$ is a module over $\tmf$, we can build a cosimplicial
resolution of $\tmf\wedge B\Sigma_3$ by $H$-modules in the category of $\tmf$-module
spectra. This greatly simplifies our computations, as the role of
the dual Steenrod algebra is played by the Hopf algebra
\[\sA:=\pi_\ast(H \wedge_{\tmf} H).\]

\begin{thm}[Henriques-Hill]
As a Hopf algebra,
\[\sA=\sA(1)_\ast\otimes E(a_2),\] where $|a_2|=9$, and $\sA(1)_\ast=\mathbb
F_3[\xi_1]/\xi_1^3\otimes E(\tau_0,\tau_1)$ is dual to the
subalgebra of the Steenrod algebra generated by $\beta$ and
$\mathcal P^1$. The elements in $\sA(1)_\ast$ have their usual
coproducts, and
\[\Delta(a_2)=1\otimes a_2+\xi_1\otimes\tau_1-\xi_1^2\otimes \tau_0+a_2\otimes 1.\]
\end{thm}

\begin{proof}
That this is a Hopf algebra follows from a slight recasting of
Adams' original analysis of the Adams spectral sequence, using the
fact that $\sA$ is flat over $H_\ast$ \cite{baker-2001-1}. We begin
with an observation of Hopkins and Mahowald, as formulated by
Behrens \cite{behrensv2}. If we let
\[C=S^0\cup_{\alpha_1} e^4\cup_{\alpha_1} e^8,\] then smashing with $\tmf$ gives
\[\tmf\wedge C=\tmf_0(2)=\BP{2}\vee\,\Sigma^8
\BP{2}.\] The middle equality demonstrates that this is
actually an $E_\infty$-ring spectrum. If we smash $\tmf\wedge C$
with $V(1)$, then we again get a ring spectrum, since the
obstruction to $V(1)$ being a ring spectrum lies in positive
Adams-Novikov filtration \cite{behrensv2}, and the homotopy of
$\tmf\wedge C$ is concentrated in filtration zero. The
Atiyah-Hirzebruch spectral sequence allows us to compute the ring
structure on homotopy, and we see that the natural map
\[S^8\to\tmf\wedge C\wedge V(1)=k(2)\vee\Sigma^8
k(2)\] behaves like a square root of $v_2$ \cite{K(2)LocalSphere}. In
other words,
\[\pi_\ast(\tmf\wedge C\wedge V(1)\big)=\mathbb
F_3[\sqrt{v_2}].\] The cofiber of $\sqrt{v_2}$ is $H$, and we have
therefore realized $H$ as a quotient of an extended $\tmf$ module by
itself via a $\tmf$ module map.

To finish the proof, we smash this cofiber sequence with $H$ over $\tmf$, giving the cofiber sequence
\[\Sigma^8 H\wedge_{\tmf}\big(\tmf\wedge C\wedge V(1)\big)\xrightarrow{\sqrt{v_2}} H\wedge_{\tmf}\big(\tmf\wedge C\wedge V(1)\big)\to H\wedge_{\tmf} H.\]
We begin by analyzing the homotopy of the first two $\tmf$ modules in this resolution:
\[\pi_\ast\Big(H\wedge_{\tmf}\big(\tmf\wedge C\wedge V(1)\big)\Big)=H_\ast\big(C\wedge V(1);\,\Z/3\big).\]

The structure of this as a graded vector space is that of
$\sA(1)_\ast$. Since $\sA$ is a commutative Hopf algebra, the classification of
Hopf algebras over a finite field ensures both that $\sqrt{v_2}$ is
zero in homotopy and that the structure of this as an algebra is as
listed \cite{MilnorMoore}. This is immediate from considering the degrees of the elements, since odd elements must be exterior classes and the element in degree $4$ must be the generator of a truncated polynomial algebra.

Since the structure map from $\tmf$ to $H$ is a map of
$E_{\infty}$ ring spectra, $\sA$ is a module over $H_\ast H$. Moreover, this map is also a map of coalgebras over $H_\ast$. Since the unit map $S^0\to\tmf$ is a $6$-equivalence, the natural map \[H\wedge H\to H\wedge_{\tmf} H\] is a $6$-equivalence. This implies that the induced map in homotopy is a Hopf algebra isomorphism in the same range, and this gives the coproducts on the elements $\tau_0$, $\tau_1$ and $\xi$.

To determine the coproduct on $a_2$, we endow $\sA$ with a filtration such that $a_2$ is primitive in the associated graded. This filtration gives rise to a spectral sequence
\[\Ext_{Gr(\sA)}(\mathbb F_3,\mathbb F_3)\Rightarrow\Ext_{\sA}(\mathbb F_3,\mathbb F_3)\] converging to the $E_2$ term of the Adams spectral sequence which computes $\pi_\ast(\tmf)$. We shall use the known computation of $\pi_\ast(\tmf)$ to deduce differentials in this algebraic spectral sequence, and this will determine the coproduct on $a_2$.

We first filter $\sA$ by letting $\sA(1)_\ast$ have filtration $0$ and
putting $a_2$ in filtration $1$. The initial piece of data needed is
the cohomology of $\sA(1)_\ast$. An elementary computation shows that as an algebra
\[\Ext_{\sA(1)_\ast}(\mathbb F_3,\mathbb F_3)=\mathbb 
F_3[v_0,v_1^3,\beta]\otimes E(\alpha_1,\alpha_2)/\big(v_0(\alpha_1,\alpha_2),\,\alpha_1\alpha_2=v_0\beta\big).\] 
This is pictorally represented in Figure \ref{ExtA1at3}.

\begin{figure}[h]
\centering \epsfig{file=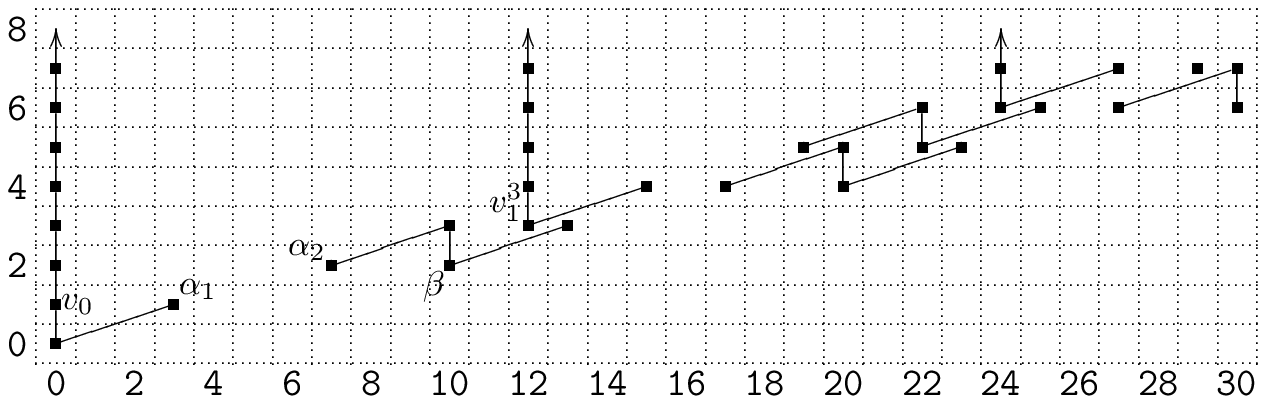,width=\textwidth}
\caption{$\Ext_{\sA(1)_\ast}(\mathbb F_3,\mathbb F_3)$} \label{ExtA1at3}
\end{figure}

Since $a_2$ is primitive in the associated graded Hopf algebra, we
know that
\[\Ext_{Gr(\sA)}(\mathbb F_3,\mathbb F_3)=\Ext_{\sA(1)_\ast}(\mathbb F_3,\mathbb F_3)[\tilde{c}_4].\]
This $\Ext$ group is the $E_1$ page of a spectral sequence
converging to the Adams $E_2$ term $\Ext_{\sA}(\mathbb F_3,\mathbb F_3)$. Since there is nothing in dimension $7$ in $\tmf_\ast$, we know that the element $\alpha_2$ must be killed. The only possible way for to achieve this is for $d_1(\tilde{c}_4)=\alpha_2$. This $E_1$ page is given
together with this necessary $d_1$ differential in Figure
\ref{tmfS0E1}.

\begin{figure}[h] \centering \epsfig{file=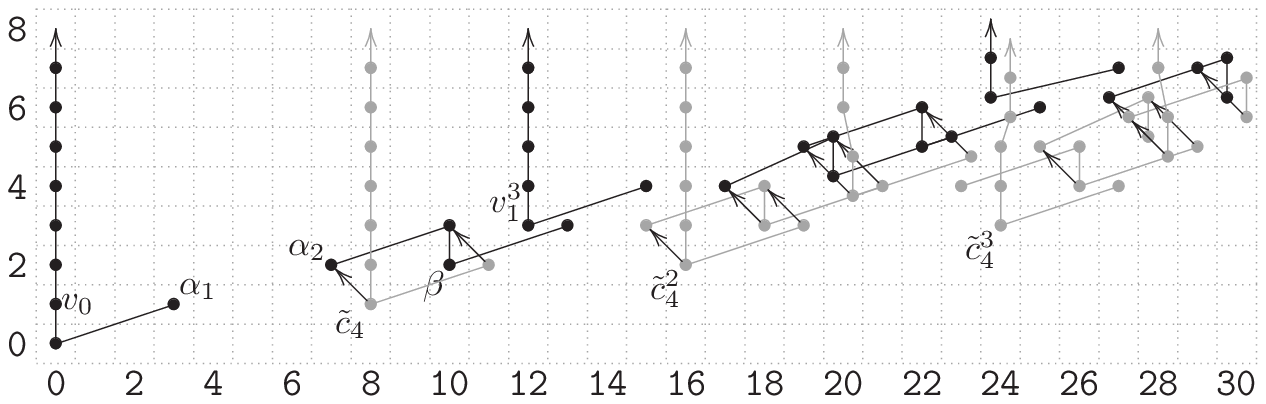, width=\textwidth}
\caption{$\Ext_{Gr(\sA)}(\mathbb F_3,\mathbb F_3)$}
\label{tmfS0E1}
\end{figure}

At this point, we rename some of the remaining elements:
\[c_4=v_0\tilde{c}_4,\quad c_6=v_1^3,\quad \Delta=\tilde{c}_4^3.\]
For completeness, we note that a similar analysis gives the $d_2$ differentials:
\[d_2([\alpha_2\tilde{c}_4^2])=v_1^3\beta,\text{ and }
d_2([v_0\tilde{c}_4^2])=v_1^3\alpha.\] The
$E_2$ page with the $d_2$ differentials is included as Figure
\ref{tmfS0E2}.

\begin{figure}[h] \centering \epsfig{file=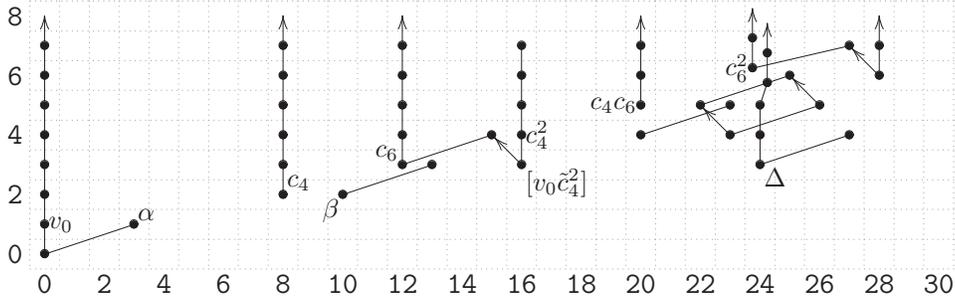, width=\textwidth}
\caption{May $E_2$ page for $\Ext_{\sA}(\mathbb F_3,\mathbb F_3)$}
\label{tmfS0E2}
\end{figure}

For the $d_1$ to have the appropriate form, we must have
\[\psi(a_2)=1\otimes a_2+a_2\otimes 1\pm(\xi_1\otimes\tau_1-\xi_1^2\otimes\tau_0).\] If the sign is negative, then we can simply replace $a_2$ by $-a_2$ to correct this.
\end{proof}

\section{Review of $\ko_\ast(\mathbb RP^\infty)$}\label{Review}
In \cite{OpsSq4}, Mahowald uses the homology of cofiber $R_2$ of the
transfer map $B\Sigma_2\to S^0$ to compute its $\ko$ homology and
the $\ko$ homology of $\mathbb RP^{\infty}$. Since the method we
will employ to handle $\tmf_\ast(B\Sigma_3)$ is similar, we quickly
review Mahowald's technique here. For this section only, all
computations will be done at the prime $2$.

\subsection{General Results and Definitions}

The homology of $R_2$ sits as an extension of the homology of
$\Sigma\mathbb RP^{\infty}$ by the homology of $S^0$, and let $e_i$ denote the generator of $H_i(R_2)$. The coaction of the dual Steenrod algebra on $H_\ast(R_2)$ is determined by the comodule structure on $H_\ast(\Sigma\mathbb RP^{\infty})$ and the coaction formula \[\psi(e_2)=\xi_1^2\otimes e_0+1\otimes e_2.\]

Let $A(1)$ be the spectrum whose cohomology is a free $\sA(1)$-module of rank $1$. Smashing $A(1)$ with $\ko$ gives a presentation of $H\Z$ as a $\ko$-module spectrum. An analysis like that of the first section reestablishes the following classical result, normally
proved using a change of rings argument.

\begin{prop}
There is a spectral sequence converging to the $\ko$ homology of a
space $X$ with $E_2$ term $\Ext_{\sA(1)_\ast}\big(\mathbb F_2,H_\ast(X)\big)$.
\end{prop}

\subsection{The $ko$ homology of $R_2$}

Mahowald's key observation was that there is a filtration of
$H_\ast(R_2)$ such that the associated graded is a sum of comodules
over $\sA(1)_\ast$ whose $\Ext$ groups are easy to compute.

\begin{prop}
There is a filtration of $H_\ast(R_2)$ such that the associated graded
is
\[Gr=Gr\big(H_\ast(R_2)\big)=\bigoplus_{k=0}^{\infty}\Sigma^{4k}h\zZ,\] where
$h\zZ$ is the $\sA(1)_\ast$ comodule dual to $\sA(1)//\!\sA(0)$.
\end{prop}

\begin{figure}[h]
\centering \epsfig{file=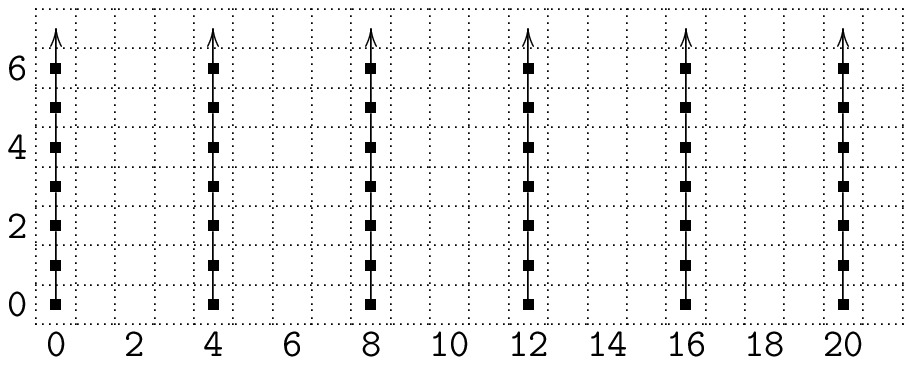}
\caption{$\Ext_{\sA(1)_\ast}(\mathbb F_2,Gr)$} \label{Gr2koRPinf}
\end{figure}

The proposition shows that if we compute $\Ext$ of $Gr$, then we see
that it is torsion free, with a $\,\zZ$ in dimensions congruent to
$0\mod 4$ (Figure \ref{Gr2koRPinf}).

Since this is concentrated in even degrees, both the algebraic
extension spectral sequence and the Adams spectral sequence
collapse. There are non-trivial extensions, though, as a
$\ko_\ast$-module.

\begin{lemma}
As a module over $\ko_\ast$,
\[ko_\ast(R_2)=\zZ_2\left[\frac{v_1^2}{4}\right].\]
\end{lemma}

\begin{proof}
An elementary cobar computation for $\Ext_{\sA(1)_\ast}(\mathbb F_2,\mathbb F_2)$ shows that the generator of the
$\zZ$ in dimension $4$ in $\ko_\ast$ is represented by
\[
[2v_1^2]=\xi_1\otimes\xi_2\otimes\xi_2+\dots.\]
If we look in the cobar complex computing $\Ext_{\sA(1)_\ast}(\mathbb
F_2,h\zZ)$, then we see that there is a class $x_7$ such that
\[x_7=\xi_1\otimes\xi_1\otimes e_5+\dots,\quad\text{ and }\quad
d(x_7)=[2v_1^2],\] where $e_5$ denotes the $5$ dimensional class in
$h\zZ$. In $H_\ast(R_2)$, the coproduct on $e_5$ is
\[\psi(e_5)=\big(\xi_1^2\otimes e_3+\xi_2\otimes
e_2+\xi_1^2\xi_2\otimes e_0\big)+\xi_1\otimes e_4+1\otimes e_5.\] This shows that
$\xi_1\otimes\xi_1\otimes\xi_1\otimes e_4$ is cohomologous to
$[2v_1^2]\otimes e_0$. The homogeneity of the homology of $R_2$ then
implies the result.
\end{proof}

\begin{remark}
This lemma shows that Mahowald and Davis' result in \cite{koTate} that $ko\wedge R_2$ splits
as a wedge of copies of $H\Z$ is not true in the category of
$ko$-module spectra.
\end{remark}

\subsection{Computing $\ko_\ast(\mathbb RP^{\infty})$}

\begin{figure}[htb]
\centering \epsfig{file=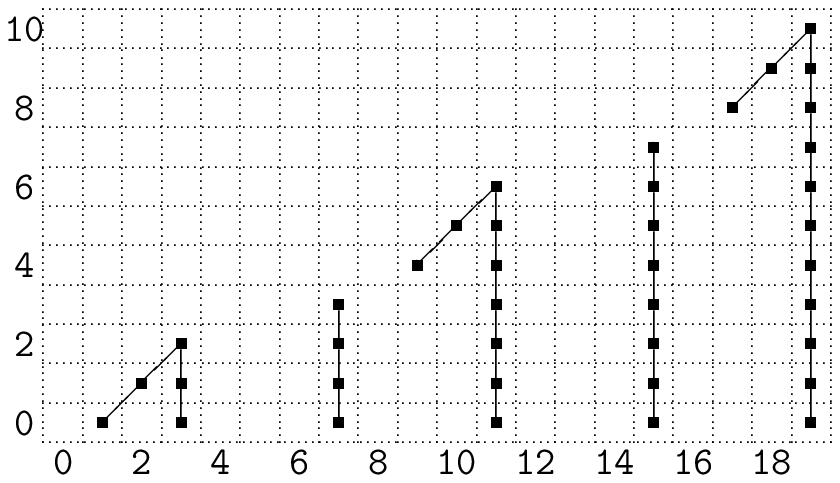} \caption{$\ko_\ast(\mathbb
RP^{\infty})$} \label{koRPinf}
\end{figure}

Finishing the argument requires looking at the long exact sequence
in $\ko$ homology for the cofiber sequence 
\[S^{0} \to R_2 \to \Sigma\mathbb RP^{\infty}.\] 
The first map is the inclusion of the
zero cell, and takes $1$ to $1$. From this, the result is easily
determined (Figure \ref{koRPinf}).

\section{The $\tmf$ Homology of the Cofiber of the Transfer $B\Sigma_3\to S^0$}\label{transfer}
Homologically, the situation at the prime $3$ is analogous to the computation at $2$. Let $R_3$ denote the cofiber of the transfer map $B\Sigma_3\to S^0$. The homology of $R_3$ sits as an extension of the homology of $\Sigma B\Sigma_3$ by the homology of $S^0$, and again let $e_i$ denote the generator of $H_i(R_3)$. The coaction of the dual Steenrod algebra on $H_\ast(R_3)$ is determined by the comodule structure on $H_\ast(\Sigma B\Sigma_3)$ and the coaction formula \[\psi(e_4)=-\xi_1\otimes e_0+1\otimes e_4.\]

The $\tmf$ analogue $h\zZ$ is again the comodule
dual to the quotient module of $\sA(1)$ by $\sA(0)$, and
the coproduct is the one induced by this structure.

\begin{lemma}\label{HR3}
$H_\ast(R_3)$ admits a filtration for which the
associated graded is
\[Gr\big(H_\ast(R_3)\big)=\bigoplus_{k=0}^{\infty}\Sigma^{12k}h\zZ.\]
\end{lemma}

\begin{proof}
In fact, this lemma is quite easy to show. The $-k^{\text{th}}$
stage of the filtration is given by taking the subcomodule generated
by the classes in dimensions $12n+1$ for all $n>k$. An elementary
computation in the cohomology of the symmetric group shows that the
associated graded is exactly what is claimed.
\end{proof}

\begin{lemma}
\[\Ext_{\sA}(\mathbb F_3,h\zZ)=\mathbb F_3\left[v_0,\frac{c_4}{3}\right].\]
\end{lemma}

\begin{proof}
To prove this lemma we apply a long sequence of spectral
sequences. First filter $\sA$ as before by letting $\sA(1)_\ast$ have
filtration $0$ and $a_2$ have filtration $1$. This filtration
extends to a filtration of $h\zZ$ in an obvious way, and we have a
spectral sequence
\[\Ext_{Gr(\sA)}(\mathbb F_3,h\zZ)\Rightarrow \Ext_{\sA}(\mathbb F_3,h\zZ).\]
As a Hopf algebra, $Gr(\sA)$ is very simple: the algebra structure
stays the same, and now $a_2$ is primitive. Now we can use the two
short exact sequences of Hopf algebras
\[ \sA(1)_\ast \to \sA\to E(a_2)\quad \text{ and }\quad E(a_2) \to \sA\to \sA(1)_\ast\]
to get a spectral sequence that converges to this $\Ext$ group and
starts with \[\Ext_{E(a_2)}\big(\mathbb F_3,\Ext_{\sA(1)_\ast}(\mathbb
F_3,h\zZ)\big).\] A final change of rings argument shows that
\[\Ext_{\sA(1)_\ast}(\mathbb F_3,h\zZ)=\mathbb F_3[v_0],\]
and this forces the result in question, since the target of any
differential on the polynomial generator is zero for degree
reasons. Again, from the previous computation of the $E_2$ term
for the Adams spectral sequence for $\tmf_\ast$ in the category of $\tmf$-modules, we see that the polynomial generator coming from $a_2$ is
$\tfrac{c_4}{3}$.
\end{proof}

\begin{figure}
\centering \epsfig{file=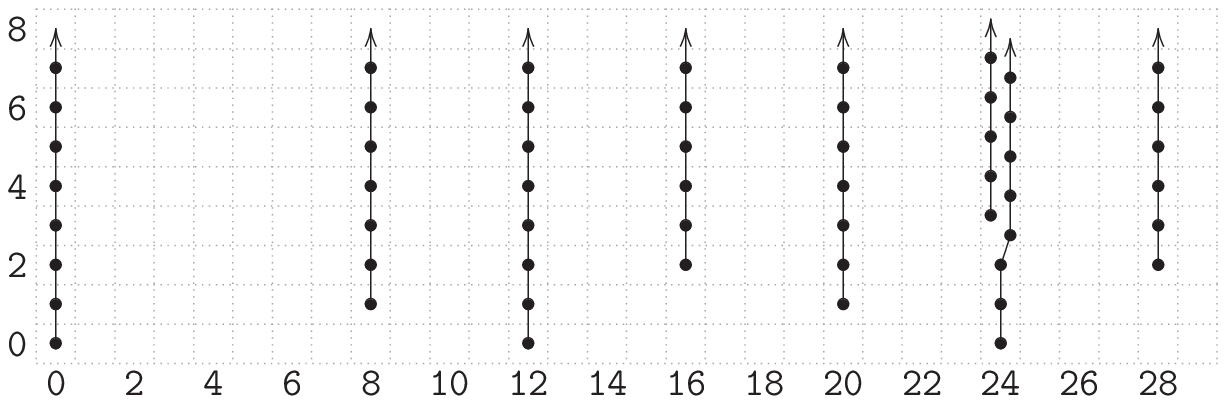}
\caption{$\Ext_{\sA}\big(\mathbb F_3,H_\ast(R_3)\big)$}
\label{Gr1tmfBSigma3-1}
\end{figure}

Since this algebra is concentrated in even degrees and since each
of the graded pieces starts an even number of steps apart, the
spectral sequence starting with $\Ext$ of the associated graded
for $H_\ast(R_3)$ collapses. We are left with the
following terms (Figure \ref{Gr1tmfBSigma3-1}).

\begin{lemma}
\[\Ext_{\sA}\big(\mathbb F_3,H_\ast(R_3)\big)=\bigoplus_{k=0}^{\infty} \Sigma^{12k}\mathbb F_3\left[v_0,\frac{c_4}{3}\right].\]
\end{lemma}

While there are no possible differentials in the Adams spectral
sequence, there are non-trivial extensions in this, viewed as a
module over $\tmf_\ast$.

\begin{thm}\label{ExtensionThm}
The Adams spectral sequence for the $\tmf$ homology of $R_3$
collapses, and as a $\tmf_\ast$-module,
\[\tmf_\ast(R_3)=\mathbb Z_3\left[\frac{c_4}{3},\frac{c_6}{27}\right].\]
\end{thm}

\begin{proof}
We show this by returning to the cobar complex. Since the homology
of $R_3$ has the very simple pattern of copies of $h\zZ$ connected
by a $\tau_0$ comultiplication on the top class in each hitting the
bottom class in the next, it will suffice to show that in the first
copy, $c_6$ on the $0$ cell is cohomologous to $27$ on the $12$
cell.

For simplicity, we will let $i_n$ denote the class in dimension $n$
in $h\zZ$. The cobar complex for
$\Ext_{\sA(1)_\ast}(\mathbb F_3,h\zZ)$ shows that there is an element
$x_{16}$ such that
\[x_{16}=\tau_0\otimes\tau_0\otimes i_{13}+\dots\text{ and } d(x_{16})=c_6\otimes i_{0}.\]
This bounding cycle can be readily found by considering
the $\Ext$ implications of the short exact sequence of comodules:
\[\mathbb F_3\{i_{0},i_{4}, i_{8}\} \to h\zZ \to \mathbb F_3\{i_5,
i_9, i_{13}\}.\]

When we add in the next copy of $h\zZ$, we change the coproduct on
$i_{13}$ to
\[\psi(i_{13})=\big(\xi_1 \otimes i_9 + \xi_1^2 \otimes i_5 + \tau_1\otimes i_8 + \xi_1\tau_1\otimes i_4 + \xi_1^2\tau_1\otimes i_{0}+1\otimes i_{13}\big)+\tau_0\otimes i_{12}.\]\
This is the only change to the coproducts in our comodule, so when
we consider again $x_{16}$ and take its boundary, the only change is
the addition of terms coming from this new term in the coproduct.
However, the only instance of $i_{13}$ in $x_{16}$ is the one coming
from $\tau_0\otimes\tau_0\otimes i_{13}$, so the real boundary is
\[ d(x_{16})=c_6\otimes i_{0}+\tau_0\otimes\tau_0\otimes\tau_0\otimes i_{12}.\] In other
words, $c_6$ on the base class is (up to a sign) $27$ times the
class in dimension $12$.
\end{proof}

\section{The $\tmf$ Homology of $B\Sigma_3$}\label{BSig3}

The most difficult of the computations now behind us, we can compute
the $\tmf$ homology of $B\Sigma_3$ by simply considering the long
exact sequence induced by applying $\tmf_\ast$ to the cofiber sequence
\[S^0\to R_3\to \Sigma B\Sigma_3.\]

The first map is the inclusion of the zero cell into $R_3$, and so
this map in $\tmf$-homology just takes $1$ to $1$. Since this is a map of
$\tmf_\ast$-modules, we see immediately that this map is injective on
elements of Adams-Novikov filtration $0$, with image

\[\mathbb Z_3\big[c_4,c_6,[3\Delta], [3\Delta^2], [c_4\Delta], [c_4\Delta^2], [c_6\Delta], [c_6\Delta^2], \Delta^3\big]/(27\Delta=c_4^3-c_6^2)\subset\mathbb Z_3[\frac{c_4}{3},\frac{c_6}{27}].\]

Additionally, since $\alpha$ and $\beta$ act as zero on all of the
classes in $\tmf_\ast(R_3)$, the kernel of this first map is the
submodule of $\tmf_\ast$ generated by $\alpha$, $\beta$ and their $\Delta$ translates. These
together establish the following theorem about the $\tmf$ homology
of $\Sigma B\Sigma_3$.

\begin{thm}
The $\tmf$ homology $\Sigma B\Sigma_3$ sits in a short exact sequence

\[0\to G_n\to \tmf_n(\Sigma B\Sigma_3) \to \widehat{\tmf}_{n-1}\to 0,\] where $\widehat{\tmf}_{n-1}$ is the subgroup of $\tmf_{n-1}$ of Adams Novikov filtration at least $1$ and $G_n$, the cofiber of the map $\tmf_n\to\tmf_n(R_3)$, is given by

\[G_{24k+12j+8i}=\begin{cases}
\mathbb Z/3\oplus\bigoplus_{m=1}^k \Z/3^{6m} & k\equiv 1,2\!\!\!\! \mod 3, i+j=0 \\
\bigoplus_{m=0}^k \Z/3^{6m+3j+i} & k\equiv 0\!\!\!\!\mod 3 \\
\bigoplus_{m=0}^k \Z/3^{6m+3j+i} & k\equiv 1,2\!\!\!\! \mod 3, i+j>0 \\
0 & \text{otherwise}
\end{cases},\] where $j<2$, and $i<3$. The sequence is split as a sequence of groups. There is a hidden $\alpha$ extension originating on the copy of $\beta^2$ in $\widehat{\tmf}_{20}$ and hitting the $\mathbb Z/3$ summand of $G_{24}$.
\end{thm}

\begin{proof}

This short exact sequence is just a restatement of the earlier
comments about the long exact sequence in $\tmf$ homology. It is
split because the elements coming from $G_n$ have Adams-Novikov
filtration $0$, and the convergence of the Adams-Novikov spectral
sequence ensures a map of groups from $\tmf_\ast(\Sigma B\Sigma_3)$ to $G_n$
which is a left inverse to this inclusion.

The structure of the groups $G_n$ is easy to show. A basis for
$\tmf_\ast(R_3)$ is given by the collection of monomials of the form
$\Delta^k\tilde{c}_6^j\tilde{c}_4^i,$ where $i<3$, and
$27\tilde{c}_6=c_6$, $3\tilde{c}_4=c_4$. This is simply because if
we can solve the relation on $\Delta$ in $\tmf_\ast(R_3)$. A basis for
the Adams-Novikov filtration $0$ subring of $tmf_\ast$ is given by the
monomials
\[\Delta^kc_6^jc_4^i \text{ for } k\equiv 0\!\!\!\mod 3\text{ or }k\equiv 1,2\!\!\!\mod 3, i+j>0,\quad [3\Delta]\Delta^k, \text{ and } [3\Delta^2]\Delta^k.\]

Recalling that
\[\Delta^kc_6^jc_4^i=3^{3j+i}\Delta^k\tilde{c}_6^j\tilde{c}_4^i\]
and collecting all terms of the same degree yields $G_n$.

The hidden extension can most readily been seen by considering the long exact sequence in $\Ext$ induced by the cofiber sequence. In this situation, $\Delta$ from the ground sphere kills $\Delta$ in the Adams $E_2$ term for $\tmf_\ast(R_3)$, and $\alpha\beta^2$ on the ground sphere survives.
\end{proof}

\begin{remark}
The proof of this theorem also shows that the transfer induces a
bijection between the elements of higher Adams-Novikov filtration
elements of $\tmf_\ast$ and the elements of $\tmf_\ast(B\Sigma_3)$ of
Adams-Novikov filtration at least one (together with the $\mathbb
Z/3$ coming from the $3$-cell). This exactly repeats the situation
at at the prime $2$, where the transfer again mapped the higher
Adams-Novikov elements in $\ko_\ast(\mathbb RP^{\infty})$ bijectively
onto those in $\ko_\ast$.
\end{remark}

\section{The $\tmf$ Homology of the Finite Skeleta of $R_3$ and $B\Sigma_3$}\label{Truncated}
For completeness, we include the $tmf$-homology of the finite skeleta of $R_3$ and $B\Sigma_3$. These computations serve as starting points for the program of Minami to detect the $3$-primary $\eta_j$ family \cite{MinamiEtaj}.

\subsection{The Skeleta of $R_3$}

Let $n=12k+i$, for $0<i\leq 12$. We wish to compute the $\tmf$-homology of $R_3^{[n]}$.

\begin{lemma}
There is a filtration of $H_\ast(R_3^{[12k+i]})$ such that the associated graded is
\[Gr\big(H_\ast(R_3^{[12k+i]})\big)=\left(\bigoplus_{n=0}^{k-1}\Sigma^{12n}h\mathbb Z\right)\oplus \Sigma^k M_{i},\] where $M_i$ is the subcomodule of $h\mathbb Z$ generated by all classes of degree at most $i$ for $i<12$, and $M_{12}$ is $M_9$ plus a primitive class in dimension $12$.
\end{lemma}

\begin{proof}
The required filtration is just the restriction of the filtration used in the proof of Lemma \ref{HR3} to the subcomodule $H_\ast(R_3^{[12k+i]})$. 
\end{proof}

The comodules $M_i$ are the homology of $R_3^{[i]}$, and this splitting result and the follow theorem demonstrates that knowing their $\tmf$-homology gives that of all finite skeleta. The proof of Theorem \ref{ExtensionThm} shows the following

\begin{thm}
As a module over $\tmf_\ast$, 
\[\tmf_\ast(R_3^{[12k+i]})=\mathbb Z_3\left[\frac{c_4}{3}\right]\{e_0,e_{12},\dots,e_{12(k-1)}\}\oplus \widetilde{M}_{i}e_{12k}/(c_6 e_{12j}-27e_{12(j+1)}),\] where $\widetilde{M}_i$ is the $\tmf$-homology of spectrum $R_3^{[i]}$.
\end{thm}

The remainder of the section will be spent computing the modules $\widetilde{M}_i$. To save space, in what follows we use two indices: $\delta$ which ranges from $0$ to $2$ and $\epsilon$ which ranges from $0$ to $1$. When these appear, it means that all possible values of the index are actually present.

\begin{prop}
The spectra $R_3^{[1]}$, $R_3^{[2]}$, and $R_3^{[3]}$ are simply $S^0$. This implies that
\[\widetilde{M}_i=\tmf_{\ast},\quad 1\leq i\leq 3.\]
\end{prop}

\begin{lemma}
The spectrum $R_3^{[4]}$ is the cofiber of $\alpha_1$. The $\tmf$-homology of this is the extension of the module generated by $[\Delta^\epsilon e_0]$ and $[\alpha e_4]$
and subject to the relations
\[ \alpha[\alpha e_4]=\beta e_0,\, \alpha[\Delta e_0]=\beta^2[\alpha e_4],\, \alpha e_0=\beta^3[\Delta^\epsilon e_0]=I[\alpha e_4]=\beta^4[\alpha e_4]\]
by the module 
\[\mathbb Z_3[c_4,c_6,\Delta]\{[3e_4],[c_4e_4],[c_6e_4].\]
The extension is determined by the two relations
\[c_4[3e_4]=3[c_4e_4]\pm c_6 e_0,\quad c_6[3e_4]=3[c_6e_4]\pm c_4^2 e_0.\] 
\end{lemma}

\begin{proof}
The Adams $E_2$ term can be readily computed to be the extension of
\[\mathbb F_3[v_0,c_4,c_6,\Delta,\beta]\{e_0\}\] by 
\[\mathbb F_3[v_0,c_4,c_6,\Delta]\{[v_0e_4],[c_4e_4],[c_6e_4]\}\oplus \mathbb F_3[\Delta,\beta]\{[\alpha e_4]\},\] subject to the relations
\[c_4[3e_4]=3[c_4e_4]\pm c_6 e_0,\quad c_6[3e_4]=3[c_6e_4]\pm c_4^2 e_0,\quad \alpha[\alpha e_4]=\beta.\] This Adams spectral sequence is a spectral module over the Adams spectral sequence for the $\tmf$-homology of the sphere, and the two differentials in the Adams spectral sequence for the sphere,
\[d_2(\Delta)=\alpha\beta^2,\quad d_3([\alpha\Delta^2])=\beta^5,\] imply that $\Delta e_0$ and $\Delta^2 e_0$ are $d_2$ cycles and that the following differentials hold:
\[d_2(\Delta[\alpha e_4])=\beta^3 e_0,\quad d_3(\alpha\Delta^2[\alpha e_4])=\beta^5[\alpha e_4].\] This last $d_3$ implies that in fact,
\[d_3(\Delta^2 e_0)=\beta^4[\alpha e_4],\] using the relation involving $\alpha$ multiplication on $[\alpha e_4]$.
\end{proof}

\begin{lemma}
The spectra $R_3^{[5]}$, $R_3^{[6]}$, and $R_3^{[7]}$ are the cofiber of the extension of $\alpha$ over the mod $3$ Moore spectrum. The $\tmf$-homology of these spectra, $\widetilde{M}_i$ is the $\tmf_\ast$ module generated by
\[[\tfrac{c_4}{3}\Delta^\delta e_0],\, [\tfrac{c_6}{3}\Delta^\delta e_0],\, [\Delta^\epsilon e_0],\, [\alpha e_4],\, [\beta e_5],\] and subject to the relations
\begin{multline*}
\alpha[\beta e_5]=\beta[\tfrac{c_4}{3}e_0],\, \alpha[\alpha e_4]=\beta e_0,\, \alpha[\Delta e_0]=\beta^2[\alpha e_4],\\ (\alpha,\beta^3) e_0=I([\alpha e_4],[\beta e_5])=\beta^4[\alpha e_4]=0.
\end{multline*}
\end{lemma}

\begin{proof}
In the long exact sequence in $\Ext$ induced by the inclusion of the $4$-skeleton into $R_3^{[5]}$, the inclusion of the $5$-cell kills the element $[v_0e_4]$. The elements $[c_4 e_4]$ and $[c_6 e_4]$ survive, and the relations in the $\Ext$ term for the $4$-skeleton ensure that in the Adams $E_2$ term for $\widetilde{M}_5$, \[v_0[c_4 e_4]=c_6 e_0,\quad v_0[c_6 e_4]=c_4^2 e_0.\] Moreover, since $\alpha$ and $\beta$ multiplications on the class $[v_0 e_4]$ are trivial, the classes $[\alpha e_5]$ and $[\beta e_5]$ survive to the Adams $E_2$ page. An elementary computation in the bar complex establishes that \[v_0[\alpha e_5]=c_4 e_0.\] 

This shows that the Adams $E_2$ page, as a module over that for $\tmf_\ast$, is
\begin{multline*}
\mathbb F_3[v_0,c_4,c_6,\Delta,\beta]\{e_0, [\tfrac{c_4}{v_0}e_0], [\tfrac{c_6}{v_0}e_0], [\alpha e_4], [\beta e_5]\}\\ /\big(\alpha[\alpha e_4]-\beta e_0, \beta[\tfrac{c_4}{v_0}e_0]-\alpha[\beta e_5], \alpha e_0,I([\beta e_5],[\alpha e_4])\big)
\end{multline*}

The differentials again follow from those in the Adams spectral sequence of $\tmf_\ast$.
\end{proof}

At this point, the patterns of extensions and differentials repeats. This makes the final computations substantially easier.

\begin{lemma}
The spectrum $R_3^{[8]}$ is the spectrum $C$ from \S \ref{HopfAlg}, where the middle cell is replaced by the mod $3$ Moore spectrum. The module $\widetilde{M}_8$ sits in a short exact sequence
\begin{multline*}
0\to \tmf_\ast\{[\tfrac{c_4}{3}\Delta^\delta e_0], [\tfrac{c_6}{3}\Delta^\delta e_0], [\Delta^\delta e_0], [\beta e_5]\}/\big((\alpha,\beta)\big([\tfrac{c_4}{3}^\epsilon\Delta^\delta e_0], [\tfrac{c_6}{3}\Delta^\delta e_0]\big),I[\beta e_5]\big)\\ \to \widetilde{M}_8\to \mathbb Z_3[c_4,c_6,\Delta]\{[3e_8],[c_4e_8],[c_6e_8]\}\to 0,
\end{multline*}
where the extension is determined by the two relations
\[c_4[3e_8]=3[c_4e_8]\pm c_4[\tfrac{c_4}{3} e_0],\quad c_6[3e_8]=3[c_6e_4]\pm c_4[\tfrac{c_6}{3} e_0].\]
\end{lemma}

\begin{proof}
The long exact sequence in $\Ext$ coming from the short exact sequence in homology induced by the inclusion of $R_3^{[5]}$ into $R_3^{[8]}$ is determined by the connecting homomophism which takes $e_8$ to $[\alpha e_4]$. The linearity of this map shows that the Adams $E_2$ term for $\widetilde{M}_8$ is an extension of
\[\mathbb F_3[v_0,c_4,c_6,\Delta]\{e_0,[\tfrac{c_4}{3} e_0], [\tfrac{c_6}{3} e_0]\}\oplus\mathbb F_3[\Delta,\beta]\otimes E(\alpha)\{[\beta e_5]\}\] by 
\[\mathbb F_3[v_0,c_4,c_6,\Delta]\{[v_0e_8],[c_4 e_8], [c_6 e_8]\},\] subject to the extensions
\[c_4[v_0 e_8]=v_0[c_4 e_8]\pm \tfrac{c_4^2}{3} e_0,\quad c_6[v_0 e_8]=v_0[c_6 e_8]\pm \tfrac{c_4c_6}{3} e_0.\]

The differentials are again determined by those of $\tmf_\ast$. The only classes which support non-trivial $\alpha$ multiplication are multiples of $[\beta e_5]$, and here, the differentials are the same as for $\widetilde{M}_5$:
\[d_2(\Delta^i[\beta e_5])=i\alpha\beta^2\Delta^{i-1}[\beta e_5],\quad d_3([\alpha\Delta^2][\beta e_5])=\beta^5[\beta e_5].\]

\end{proof}

\begin{lemma}
The spectra $R_3^{[9]}$, $R_3^{[10]}$, and $R_3^{[11]}$ are the cofiber of the map from $\Sigma^4 C(\alpha)$ to $C$ which is multiplication by $3$ on the $4$ and $8$ cells. The module $\widetilde{M}_9$ can be expressed via the short exact sequence
\[0\to\tmf_*\{[\alpha e_9]\}\to\widetilde{M}_9\to\mathbb Z_3\left[\frac{c_4}{3}\right]e_0\to 0,\] where the only extension is given by \[c_6 e_0=9[\alpha e_9].\]
\end{lemma}

\begin{proof}
The cofiber sequence coming from the inclusion of $R_3^{[8]}$ into $R_3^{[9]}$ induces a long exact sequence on $\Ext$. The connecting homomorphism is
\[e_9\mapsto [v_0 e_8]+[\tfrac{c_4}{3}e_0].\] This is a map of modules over the Adams $E_2$ term for $\tmf_\ast$, and just as before, the element $[\alpha e_9]$ is in the kernel of this map. This gives hidden extensions analogous to the ones for $\widetilde{M}_4$ and $\widetilde{M}_5$ in the Adams $E_2$ term for $\widetilde{M}_9$:
\[\alpha[\alpha e_9]=\beta e_5,\quad v_0[\alpha e_9]=[\tfrac{c_6}{3}e_0].\] The $c_4$ and $c_6$ extensions coming from $[v_0 e_8]$ give two more extensions:
\[v_0[c_4 e_8]=c_4[\tfrac{c_4}{3}e_0],\quad v_0[c_6 e_8]=c_4[\tfrac{c_6}{3} e_0].\] This establishes that the Adams $E_2$ term is given by the extension of
\[\mathbb F_3[v_0,c_4,c_6,\Delta]\{[\alpha e_9]\}\] by
\[\mathbb F_3[v_0,c_4,c_6,\Delta]\{e_0,[\tfrac{c_4}{v_0} e_0],[\tfrac{c_4^2}{v_0^2} e_0]\},\]
where $c_6 e_0=v_0^2 [\alpha e_9].$
Just as before, the ordinary Adams differentials determine the differentials, recalling that $[\tfrac{c_6}{v_0^2} e_0]=[\alpha e_9]$:
\[d_2(\Delta^k[\tfrac{c_6}{v_0^2} e_0])=k\alpha\beta^2\Delta^{k-1}[\tfrac{c_6}{v_0^2} e_0]=\beta^2[\beta e_5], d_3(\Delta^2[\beta e_5])=\beta^5[\alpha e_9].\] The Adams differentials here preserve the exact sequence, and this establishs the statement of the Lemma.
\end{proof}

\begin{remark}
For completeness, we note that if we were to include a $13$-cell, attaching it to the $9$-cell via $\alpha$, then the attaching map in long exact sequence in $\tmf$ homology would take the copy of $\tmf_\ast$ coming from the $13$-cell isomorphically onto the factor $\tmf_\ast\{[\alpha e_9]\}.$
\end{remark}

\begin{prop}
Since the twelve dimensional class is primitive in $M_{12}$, we conclude that as a $\tmf_\ast$-module, 
\[\widetilde{M}_{12}=\widetilde{M}_9\oplus \Sigma^{12}\tmf_\ast.\]
\end{prop}

\subsection{The Skeleta of $B\Sigma_3$}

The analysis of the preceding section allows us to completely determine the structure of the groups $\tmf_\ast(B\Sigma_3^{[n]})$. However, due to the complexity of the combinatorial problem, explicit demonstration of these groups in unenlightening. We instead present the following theorem concening bounds on the orders of these groups.

\begin{thm}
If $n=12k+i$, then $3^{3k+2}$ annihilates the torsion subgroup of $\tmf_\ast(B\Sigma_3^{[n]})$. Moreover, if $i\geq 5$, then there are elements of order exactly $3^{3k+1}$, and if $i\geq 9$, then there are elements of order exactly $3^{3k+2}$.
\end{thm}

\begin{proof}
This is immediate with the consideration that the large torsion subgroups are generated by high powers of $\tfrac{c_6}{27}$. If we consider only a finite skeleton of $B\Sigma_3$, then we include only finitely many powers of this element. The largest such element occurs in dimension $12k$. If $i$ is at least 5, then we have the element $\tfrac{c_4}{3}$ on this element. If $i$ is at least $9$, then we have the element $\tfrac{c_4^2}{9}$ on this element. These provide the elements of exact order.
\end{proof}

\section{The $\Sigma_3$ Tate Homology of $\tmf$}\label{tate}
The analysis used to compute the $\tmf$ homology of $R_3$ applies to
compute the homotopy of
\[\tmf^{t\Sigma_3}=\left(\tmf\wedge B\Sigma_3\right)_{-\infty}=\lim_{\longleftarrow}\big(\tmf\wedge {(B\Sigma_3)}_{-n}\big).\]

\subsection{Computation of the Homotopy}
A mod $3$ form of James periodicity shows that as $\sA(1)_\ast$-comodules,
\[H_\ast\big({(B\Sigma_3)}_{-12k+3}\big)=\Sigma^{-12k}H_\ast\big({(B\Sigma_{3})}_{3}\big).\]
The Adams spectral sequence argument in \S \ref{BSig3} shows that the map
\[\pi_\ast\big(\tmf\wedge (B\Sigma_3)_{-12(k+1)+3}\big)\to\pi_\ast\big(\tmf\wedge (B\Sigma_3)_{-12k+3}\big)\] is surjective on the $G_\ast$ summand and zero on the $\widehat{\tmf}_\ast$ summand. This implies that there are no $\lim^1$ terms coming from the inverse system of homotopy groups. Moreover, this is a system of $\tmf_\ast$-modules, and considering the action of $c_4$ and $c_6$ in each of the modules in the inverse system allows us to conclude

\begin{thm}
The homotopy of the $\Sigma_3$ Tate spectrum of $\tmf$ is an
indecomposable $\tmf_\ast$ module, and
\[\pi_\ast(\tmf^{t\Sigma_3})=\Sigma^{-1}\mathbb Z_3\left[\frac{c_4}{3},\left(\frac{c_6}{27}\right)^{\pm 1}\right].\]
\end{thm}

\subsection{A Conjectural Non-splitting Result}
We wish to establish a limit argument using the Adams spectral sequence to show that this spectrum does not split. When we consider the effect of homology on the James periodicity result, then it shows that there is a filtration of $H_\ast({B\Sigma_3})_{-\infty}$ such that the associated graded is
\[Gr\big(H_\ast({B\Sigma_3})_{-\infty}\big)=\bigoplus_{k=-\infty}^{\infty}\Sigma^{12k-1}h\zZ.\]
This implies that the limit Adams spectral sequence for the homotopy of
$\tmf^{t\Sigma_3}$ collapses, reaffirming the previous result. However, analysis of the Adams $E_2$ term shows that the Adams filtrations seem to be wrong for a splitting result analogous to Mahowald and Davis' result. The classes in dimensions $3$ mod $12$ all have Adams
filtration at least $2$, whereas $v_1$ should lie in Adams filtration $1$.

\begin{conjecture}
There does not exist a splitting of the form
\[\tmf^{t\Sigma_3}=\bigvee_{k=-\infty}^{\infty}\Sigma^{12k-1}\BP{1}.\]
\end{conjecture}

\section{The Conjectural Case for Higher Primes}\label{eop-1}

A similar result is conjectured to hold for the $p$-local case
with $\eo_{p-1}$, where $\eo_{p-1}$ is an $E_{\infty}$ ring spectrum which $K(p-1)$-localizes to $EO_{p-1}$ and whose homotopy groups are determined by the Gorbounov-Hopkins-Mahowald Hopf algebroid without inverting $\Delta$ or completing \cite{GorMah}. First, we should have a similar Hopf algebra.

\begin{conjecture}
As a Hopf algebra,
\[\pi_\ast(H\Z/p\wedge_{\eo_{p-1}} H\Z/p)=\sA(1)\otimes E(a_2,\dots,a_{p-1}),\]
where again $\sA(1)$ is dual to the subalgebra generated by
$\beta$ and $\mathcal P^1$, and where $|a_i|=2i(p-1)+1$. The
elements in $\sA(1)$ again have their usual coproducts, while
\[\psi(a_j)=\sum_{k=0}^{j}\tfrac{1}{k!} \xi_1^k\otimes a_{j-k}+a_j\otimes 1,\]
where $a_1=\tau_1$ and $a_0=\tau_0$.
\end{conjecture}

\begin{remark}
Since the spectra $\eo_{p-1}$ are not known to exist for $p>3$, we
can only comment that this follows the huristic pattern of $\eo_2$,
and if we invert $\Delta$, the class corresponding to
$[a_2]^p=\left(\tfrac{c_4}{p}\right)^p$ after running the
corresponding Adams spectral sequence, then we get the homotopy of
$EO_{p-1}$.
\end{remark}

\begin{proof}[Indicative Sketch]
The spectrum $EO_{p-1}$ is the homotopy fixed points of $E_{p-1}$ under an action of an extension of $\mathbb Z/p$ by $\mathbb Z/{(p-1)^2}$. Since $E_{p-1}$ is a $p$ complete spectrum, the prime to $p$ part of the group serves only to carve out an ``Adams summand'' for $EO_{p-1}$. The $p$-cell spectrum \[C=S^0\cup_{\alpha_1}e^{2(p-1)}\cup\dots\cup_{\alpha_1}e^{2(p-1)^2},\] 
when smashed with $EO_{p-1}$ undoes the $\Z/p$ homotopy fixed points, resulting in a torsion free spectrum that is the $K(p-1)$-localization of a wedge of copies of $\BP{p-1}$. This implies that just as in the case of $p=2$ or $p=3$, $\eo_{p-1}\wedge C$ should split as a wedge of copies of $\BP{p-1}$. The determination of the number of exterior classes and their coproducts come from considering the implications in homotopy of such a splitting, just as was done in \S \ref{HopfAlg} for the prime $3$.
\end{proof}

Assuming the proposition, most of the results true for the prime $3$
hold generically. If we again consider the cofiber $R_p$ of the
transfer map $B\Sigma_p\to S^0$, then there is an analogue to Lemma \ref{HR3}

\begin{proposition}
There is a filtration of $H_\ast(R_p)$ such that the associated graded is
\[Gr\big(H_\ast(R_p)\big)=\bigoplus_{k=0}^{\infty}\Sigma^{2p(p-1)k}h\zZ.\]
\end{proposition}

The same argument that showed that $\Ext_{\sA}$ of this was
torsion free works at other primes, so we see that $\Ext_{\sA_p}\big(\mathbb F_p,H_\ast(R_p)\big)$ is an evenly generated polynomial algebra with generators corresponding to $p$ and certain fractional multiplies of rational generators of $eo_{p-1\ast}$. The extension problems can be similar solved.

We conjecture that $\eo_{p-1\ast}(R_p)$ is again indecomposable as a module over $eo_{p-1\ast}$.
We moreover conjecture that the $eo_{p-1}$ image of the transfer
map again contains all of the higher Adams-Novikov filtration
elements, since these are generated by $\alpha$ and $\beta$, and
these elements will again not be present in $eo_{p-1\ast}(R_p)$.

\bibliography{biblio-1}

\end{document}